\documentclass[12pt,twoside]{article} 
\usepackage{amsmath} 
\usepackage{amsfonts} 
\usepackage{amsthm} 
\usepackage{amssymb} 
\usepackage{graphicx} 
\usepackage{multicol}

\begin{document} 

\title{{\normalsize{\bf Quasi-ribbon but non-ribbon surface-link of trivial components}}} 
\author{{\footnotesize Akio Kawauchi}\\ 
{\footnotesize{\it Osaka Central Advanced Mathematical Institute, Osaka Metropolitan University}}\\ 
{\footnotesize{\it Sugimoto, Sumiyoshi-ku, Osaka 558-8585, Japan}}\\ 
{\footnotesize{\it kawauchi@omu.ac.jp}}} 
\date\, 

\maketitle 

\vspace{0.25in} 
\baselineskip=10pt 

\newtheorem{Theorem}{Theorem}[section] 
\newtheorem{Conjecture}[Theorem]{Conjecture} 
\newtheorem{Lemma}[Theorem]{Lemma} 
\newtheorem{Sublemma}[Theorem]{Sublemma} 
\newtheorem{Proposition}[Theorem]{Proposition} 
\newtheorem{Corollary}[Theorem]{Corollary} 
\newtheorem{Claim}[Theorem]{Claim} 
\newtheorem{Definition}[Theorem]{Definition} 
\newtheorem{Example}[Theorem]{Example}

\begin{abstract} 
An operation on a surface-link preserving the surface-knot components called a 
quasi-equivalence is introduced, which induces an equivalence relation on the 
surface-links. 
A quasi-ribbon surface-link is a surface-link quasi-equivalent to a ribbon surface-link.  
It is shown that a surface-link of trivial components is a quasi-ribbon surface link 
if it has only at most one non-sphere component, 
whereas there exist examples that are not quasi-ribbon surface-links of trivial  components when it has at least two non-sphere components.
For every  closed disconnected, orientable surface, a quasi-ribbon but non-ribbon 
surface-link of trivial components is constructed.
This construction is done by applying a positive solution to Cochran's conjecture,  
meaning that the sphere-knot obtained from the anti-parallel sphere-link of every  
non-ribbon sphere-knot by surgery along an  S-nontrivial or S-trivial fusion 
1-handle is a  non-ribbon or trivial sphere-knot, respectively.

\phantom{x}



\end{abstract} 

\baselineskip=15pt 

\phantom{x} 

\noindent{\bf 1. Introduction} 

Let ${\mathbf F}$ be a (possibly disconnected) closed oriented surface.  
An ${\mathbf F}$-link in the 4-sphere $S^4$ is the image of a smooth embedding 
${\mathbf F} \to S^4$. 
A surface-link in $S^4$ is an ${\mathbf F}$-link for some ${\mathbf F}$. 
If ${\mathbf F}$ consists of some copies of the 2-sphere $S^2$, then it is also called an 
$S^2$-link. 
When ${\mathbf F}$ is connected, they are also called an ${\mathbf F}$-knot, 
a surface-knot, and an $S^2$-knot, respectively.  
A  trivial surface-link is a surface-link $F$ in $S^4$ which bounds disjoint  
handlebodies smoothly embedded in $S^4$.  
A  1-handle  on a surface-link $F$ in $S^4$ is a 1-handle $h$ on $F$ embedded smoothly in $S^4$, which is  called a {\it self 1-handle} on $F$ if the attaching part of $h$ belongs to the same component of $F$.  Otherwise, the 1-handle $h$ is called a 
{\it fusion} 1-handle on $F$. 
Let $F(h)$ be the surface-link obtained from $F$ by surgery along a 1-handle system (i.e., a  system of disjoint 1-handles) $h$. 
A  ribbon surface-link is the surface-link $F=O(h^O)$ in $S^4$ obtained  
from a trivial $S^2$-link $O$ by surgery along a 1-handle system $h^O$ on $O$, 
\cite{[1], [2]}.  
An  O2-handle pair on a surface-link $F$ in $S^4$ is a pair  
$(D\times I, E \times I)$ of 2-handles  
$D\times I$, $E\times I$ on $F$ in $S^4$ which intersect orthogonally only with the attaching parts  $(\partial D)\times I$, $(\partial E)\times I$ to $F$  so that  
the intersection $Q=(\partial D)\times I\cap (\partial E)\times I$ is a square, \cite{[3]}. 
Let $F_0$ be a surface-knot in $S^4$, and $h_0$ a trivial 1-handle on $F_0$ (namely, 
the core arc $\alpha_0$ of $h_0$ is an interior push of an arc in $F_0$). 
Then there is 
an O2-handle pair $(D\times I, E\times I)$ on the surface-knot $F_0(h_0)$ in $S^4$. 
For any O2-handle pair $(D\times I, E\times I)$ on $F_0(h_0)$, 
the surface-knot $F_0(h_0)(D\times I, E\times I)$ is not only equivalent to the  
surface-knots $F_0(h_0)(D\times I)$ and $F_0(h_0)(E\times I)$ by 
the common 2-handle property, but also equivalent 
to the original surface-knot $F_0$, \cite{[3],[4],[5]}.
More generally, let $F$ be a surface-link  with a component $F_0$ in $S^4$, and 
$h_0$ a self 1-handle on $F$ in $S^4$ such that $h_0$ is a trivial 1-handle 
on $F_0$. Then the 1-handle $h_0$ is called a {\it self-trivial} 1-handle on $F$. 
Let $(D\times I, E\times I)$ be an O2-handle pair on $F_0(h_0)$ in $S^4$ such that the 2-handle  $D\times I$ does not meet the complementary surface-link 
$F\setminus F_0$ while the 2-handle $E\times I$ may meet $F\setminus F_0$, which is referred to as a {\it quasi-O2-handle pair} on $F$. 
The  surface-link $F'=F(h_0)(D\times I)$ is a (possibly new) 
surface-link consisting of  the components equivalent to the original  
components of $F$.  This correspondence $F\to F'$ is called an {\it 
elementary quasi-equivalence}.    
A correspondence $F\to F^*$ for a surface-link $F^*$ is called a {\it quasi-equivalence} if there is a sequence of  surface-links 
$F^{(i)}\, (i=0,1,2,\dots, n)$  with $F^{(0)}=F$ and $F^{(n)}=F^*$ such that 
the correspondence $F^{(i-1)}\to F^{(i)}$ is an elementary quasi-equivalence
for all $i>0$. 
{\it It would be interesting to observe that the inverse correspondence 
\mbox{$F^*\to F$} is also a quasi-equivalence and thus, the 
quasi-equivalence on a surface-link is an equivalence relation.} 
This is because the 
2-handle  $D\times I$ on $F(h_0)$ is regarded as 
a self-trivial 1-handle  $h'_0$ on $F'$ with $F'=F(h_0)(D\times I)$ and 
the self-trivial 1-handle $h_0$ on $F$ is regarded as a 2-handle    
$D'\times I$ on $F'(h'_0)$ constituting 
a quasi-O2-handle pair  $(D'\times I, E'\times I)$ on $F'(h'_0)$ such that 
$F=F'(h'_0)(D'\times I)$. 
A {\it quasi-ribbon surface-link} is a surface-link quasi-equivalent to 
a ribbon surface-link. 
A  ribbon surface-link is a 
quasi-ribbon surface-link and a quasi-ribbon surface-link 
consists of ribbon surface-knot components. 
A quasi-ribbon surface-knot is a ribbon surface-knot. 
The following theorem is shown.  

\phantom{x} 

\noindent{\bf Theorem~1.1.} Every surface-link $F$ of trivial components with only  at most one non-sphere component  in $S^4$  is a quasi-ribbon surface-link.

\phantom{x} 

Theorem~1.1 corrects the author's earlier statement, \cite[Theorem~1]{[6]}, 
\cite[Theorem~1.1, Corollary~1.3]{[7]}.
The following theorem is a strengthened version of \cite[Theorem~2]{[6]} and is shown in Section~3.

\phantom{x} 

\noindent{\bf Theorem~1.2.}  For every disconnected closed oriented surface 
${\mathbf F}$ with at least two non-sphere components, there is a pair of 
${\mathbf F}$-links $F$ and $F'$ in $S^4$ both of trivial components with 
the same fundamental group up to meridian-preserving isomorphisms such that 
$F$ is a ribbon surface-link but $F'$ is not a  quasi-ribbon surface-link.  

\phantom{x} 

For the proof of Theorem~1.2, a generalization of the null-homotopic Gauss sum invariant of a  surface-knot to a surface-link is discussed, \cite{[8]}.  
The following theorem is the main result of this paper. 

\phantom{x} 

\noindent{\bf Theorem~1.3.} For every disconnected closed oriented surface 
${\mathbf F}$, there is a quasi-ribbon but non-ribbon ${\mathbf F}$-link of trivial components in $S^4$.

\phantom{x} 

Theorem~1.3  generalizes Ogasa's observation, \cite{[9]}.
The proof of Theorem~1.3 is done in Section~4  by applying a positive answer to Cochran's conjecture discussed there, \cite{[10]}.

\phantom{x} 

\noindent{\bf 2. Proof of Theorem~1.1} 

A {\it semi-unknotted multi-punctured handlebody 
system}, or simply a {\it SUPH system} for a surface-link $F$  
in $S^4$ is a compact oriented 3-manifold $W$ smoothly embedded in  
$S^4$ such that $W$ is a handlebody system with a finite number of open 3-balls removed  
and the boundary $\partial W$ of $W$ is given by $\partial W = F \cup O$ for a trivial $S^2$-link $O$  in $S^4$, \cite{[4]}. 
A typical SUPH system $W$ is constructed  
from a ribbon surface-link $F$ defined from a trivial $S^2$-link $O$  
and a 1-handle system $h^O$ on $F$ as the union $O\times[0,1]\cup h^O$  
for a normal collar $O\times [0,1]$ of $O$ in $S^4$ with $O\times \{0\}=O$ 
 where $h^O$ does not meet $O\times [0,1]$ except for the attaching part to 
 $O$ and $\partial W = F \cup O \times\{1\}$. 
For a SUPH system $W$ with $\partial W = F \cup O$, there is a proper arc system 
$\alpha$ in $W$ spanning $O$ such that a regular neighborhood $N(O\cup\alpha)$ of the union $O\cup\alpha$  in $W$ is diffeomorphic to the closed complement  
$\mbox{cl}(W\setminus c(F\times[0,1]))$ of a boundary collar $c(F\times[0,1])$ of $F$ in $W$. This pair  
$(O,\alpha)$ is called a {\it chorded sphere system} of the SUPH system $W$. 
By replacing $\alpha$  
with a 1-handle system $h^O$ attaching to $O$ with core arc system $\alpha$, the surface-link $F$ 
is a ribbon surface-link defined by $O$ and $h^O$. In other words, giving a SUPH system $W$ with  
$\partial W = F \cup O$ is the same as saying that the surface-link $F$ is a ribbon 
surface-link with sphere system $O$.  
A {\it multi-fusion SUPH system} of a SUPH system $W$ with $\partial W = F \cup O$  
in $S^4$ is a SUPH system for $F$ in $S^4$ obtained from $W$ by deleting an open regular  neighborhood of a disjoint simple proper arc system in $W$ spanning $O$.  
A {\it multi-fission SUPH system} of a SUPH system $W$ with $\partial W = F \cup O$  
in $S^4$ is a SUPH system for $F$ in $S^4$ obtained from $W$ by adding a  2-handle  
system $h^2$ on $O$ disjoint from $W$ except for the attaching part in $O$ where $h^2$ is  taken in the complement  $B\setminus D_W$ for a disjoint 3-ball system $B$ bounded by $O$ in $S^4$ with $(B\setminus O)\cap W=D_W$ a disjoint 2-disk system. 
For any two SUPH systems $W$ and $W'$ for a ribbon surface-link $F$ in $S^4$, 
there is an orientation-preserving diffeomorphism $f$ of $S^4$ sending  $W$ to
a multi-fusion SUPH system $W^{**}$ of a multi-fission  SUPH system $W^*$ of 
the SUPH system  $W'$, \cite[Theorem~A.1]{[11]}.  
The following theorem, giving a characterization of quasi-ribbon surface-links, corrects  an earlier characterization of when a surface-link of ribbon components is 
a ribbon surface-link, \cite[Theorem~1.4]{[4]}.

\phantom{x}

\noindent{\bf Theorem~2.1.}  A surface-link $F$ in $S^4$ is a quasi-ribbon surface-link 
if and only if  there is a self-trivial 1-handle system $h$ on $F$ 
such that the surface-link  $F(h)$  is a ribbon surface-link  in $S^4$.

\phantom{x} 

\noindent{\it Proof.} 
A  surface-link  $F$ is a quasi-ribbon surface-link if and only if 
there is a ribbon surface-link $F'$ with the same components as $F$ such that 
$F=F'(h')(D\times I)$ for a self-trivial 1-handle system $h'$ on $F'$ 
and a quasi-O2-handle pair system $(D\times I, E\times I)$ on $F'(h')$. 
Then the 2-handle system $D\times I$ on the ribbon surface-link $F'(h')$ 
is regarded as a self-trivial 1-handle system $h$ on $F$ 
such that $F(h)=F'(h')$. This shows the \lq\lq only if part\rq\rq  of the theorem. 
Conversely, assume  that $F(h)$ is a ribbon surface-link  
consisting of the ribbon surface-knots $F_i(h_i)\, (i=1,2,\dots, r)$ 
where $F_i$ is  a surface-knot component of  $F$ and $h_i$ is
the self-trivial 1-handle system on $F_i$ with $h=\cup_{i=1}^r h_i$.  
By the ribbonness of $F_i(h_i)$ and the triviality of $h_i$ on $F_i$,  
the surface-knots  $F_i\, (i=1,2,\dots, r)$ are ribbon surface-knots, \cite{[4]}. 
Let $W$ be a SUPH system for $F(h)$ which is a union of  SUPH systems  $W_i$ for 
$F_i(h_i)\, (i=1,2,\dots, r)$. 
Let $W^0_i$ be a SUPH system for the ribbon surface-knot $F_i$. 
The trivial 1-handle system $h_i$ on $F_i$ is obtained so that the core arc system 
$\alpha_i$ of $h_i$ is obtained from an arc system $\alpha^0_i$ in $F_i$ 
by pushing  the interior of the arc system $\alpha^0_i$ outside of $W^0_i$. 
Then the union $W'_i=W^0_i \cup h_i$ is a SUPH system for $F_i(h_i)$ which admits  
an O2-handle pair system 
$(D_{h_i}\times I, E_i\times I)$ on $F_i(h_i)$ such that the core disk system 
$D_{h_i}$ of the 2-handle system $D_{h_i}\times I$ is a transverse disk system 
of $h_i$ and the core disk system  $E_i$ of the 2-handle system $E_i\times I$  does 
not meet $W'_i$.  For the SUPH systems $W_i$ and $W'_i$ for $F_i(h_i)$, 
there is an orientation-preserving diffeomorphism $f_i$ of $S^4$ sending  $W_i$ to 
a multi-fusion SUPH system $W^{**}_i$ of a multi-fission  SUPH system $W^*_i$ of 
 the SUPH system  $W'_i$. 
Then the disk system $D_{h_i}$  is in $W'_i$ and hence in $W^*_i$, which is deformed  into a disk system $D^*_i$ in $W^{**}_i$ with $\partial D^*_i=\partial D_{h_i}$. 
Since the core arc system $\alpha_i$ of $h_i$ is an interior push of the arc system 
$\alpha^0_i$ in $F_i$,  
the disk system $E_i$ is made disjoint from a disjoint 3-ball system $B_i$ bounded by 
the trivial $S^2$-link $O_i= \partial W^0_i\setminus F_i$ by sliding 
the arc system $\alpha^0_i$ along $F_i$. 
This means that the disk system $E_i$ meets $W^{**}_i$ only with the boundary loop system $\partial E_i$. 
The disk systems $D'_i= f_i^{-1}(D^*_i)$ and $E'_i= f_i^{-1}(E_i)$
constitute an O2-handle system $(D'_i\times I, E'_i\times I)$ on $F_i(h_i)$ 
such that $D'_i$ is a proper disk system in $W_i$ and $E'_i$ meets $W_i$ only with the boundary loop system $\partial E'_i$. 
Since the pair $(D'_i\times I, E'_i\times I)$ is  a quasi-O2-handle pair system on 
$F(h_i)$, 
there is a quasi-O2-handle pair system $(D'\times I, E'\times I)$ on the ribbon 
surface-link $F(h)$ 
by taking $D'=\cup_{i=1}^r D'_i$ and $E'=\cup_{i=1}^r E'_i$ such that 
$F'=F(h)(D'\times I)$ is a ribbon surface-link with a SUPH system obtained from 
the SUPH system $W$ for $F(h)$ by cutting along the proper disk system $D'$. 
Thus, $F$ is quasi-equivalent to the ribbon surface-link $F'$, and a quasi-ribbon 
surface-link. This completes the proof of Theorem~2.1.  

\phantom{x}

The following lemma is shown in \cite[Lemma~1.2]{[7]} (although similar ideas are found in the incorrect claims, \cite[Theorem~1.4]{[4]} and \cite[Theorem~1]{[6]}).  

\phantom{x} 

\noindent{\bf Lemma~2.2.} For every surface-link $F$ in $S^4$ with at most one 
non-sphere component, there is a self 1-handle system $h$ on $F$ such that the 
surface-link $F(h)$ is a ribbon surface-link in $S^4$.

\phantom{x} 

The proof of Theorem~1.1 is done as follows.

\phantom{x} 

\noindent{\it 2.3: Proof of Theorem~1.1.} 
Every self 1-handle system  on every surface-link $F$ of trivial components 
is a  self-trivial 1-handle system on $F$, \cite{[12]}.
Thus, Theorem~1.1 is obtained from Theorem~2.1 and Lemma~2.2, completing the proof of Theorem~1.1.
 
\phantom{x} 

\noindent{\bf 3. Proof of Theorem~1.2} 

The {\it quadratic function} $\eta:H_1(K;Z_2)\to Z_2$ of a surface-knot $K$ in $S^4$  
is defined as follows.  
For a loop $\ell$ on $K$, let $d$ be a compact (possibly non-orientable) surface in $S^4$  with $d\cap K=\partial d=\ell$. 
The value $\eta([\ell])$ is defined by the $Z_2$-self-intersection number 
$\mbox{Int}(d,d) \mod{2}$  
with respect to the framing of the surface $K$ which is independent of a choice of $d$ by  calculation.  
The function $\eta:H_1(K;Z_2)\to Z_2$ is a $Z_2$-quadratic function with the identity  
\[ \eta(x+y)=\eta(x)+\eta(y) + x\cdot y,\] 
where $x,y \in H_1(K;Z_2)$ and $x\cdot y$ denotes the $Z_2$-intersection number of $x$ and $y$ in $K$.  
A loop $\ell$ on $K$ is {\it spin} or {\it non-spin} according to whether $\eta([\ell])$ is $0$ or $1$, respectively.  
For a surface-link $F$ in $S^4$, the {\it quadratic function} $\eta:H_1(F;Z_2)\to Z_2$ 
of $F$  is defined to be the split sum of the quadratic functions 
$\eta_K:H_1(K;Z_2)\to Z_2$ for all the 
components $K$ of $F$. This quadratic function may be identified with  
the quadratic function $\eta_{\#}:H_1(F_{\#};Z_2)\to Z_2$ of a surface-knot $F_{\#}$  
in $S^4$ which is a fusion of $F$ along a fusion 1-handle system on $F$ under  
a canonical isomorphism $\iota:H_1(F;Z_2)\to H_1(F_{\#};Z_2)$.  
To see this, let $\ell$ be a loop in a component $K$ of $F$, $d$  
a compact surface in $S^4$ with $d\cap K=\partial d=\ell$, and $F_{\#}$ the 
connected surface-knot in $S^4$ obtained from $F$ 
by surgery along  a fusion 1-handle system  $h$ on $F$.  
Every transverse intersection point between 
$d$ and $F\setminus K$ in $S^4$ can be moved into $K\setminus h\cap K$ through $F_{\#}$, so that the compact surface $d$ is modified into a compact surface $d_K$ in $S^4$ with $d_K\cap F= \partial d_K=\ell \cup o_K$ for a trivial loop system  $o_K$ in $K$. Since every loop of $o_K$ is a spin loop in $F_{\#}$,  the identities  
\[\eta_{\#}(\iota([\ell]))=\eta_{\#}(\iota([\ell])+ \iota([o_K]))=\eta([\ell]),\] 
hold, showing the identification of $\eta$ to $\eta_{\#}$.
Let $\Delta(F;Z_2)$ be the subgroup of $H_1(F;Z_2)$ consisting of  
an element represented by a loop $\ell$ in $F$ which bounds an immersed disk $d$ in $S^4$  with $d\cap F=\ell$.  
The restriction $\xi:\Delta=\Delta(F;Z_2)\to Z_2$ of the quadratic function $\eta$ 
on $H_1(F;Z_2)$ is  
called the {\it null-homotopic quadratic function} of the surface-link $F$.  
The {\it null-homotopic Gauss sum} of $F$ is the Gauss sum $GS_0(F)$ of $\xi$ 
defined by  
\[GS_0(F)=\sum _{x\in \Delta} \exp(\xi(x)\pi\sqrt{-1}).\] 
This number $GS_0(F)$ is an invariant of a surface-link $F$, which is calculable as  shown  for the case of a surface-knot, \cite{[8]}.  
The following lemma is used for the proof of Theorem~1.2.

\phantom{x}

\noindent{\bf Lemma~3.1.}
If $F$ is a quasi-ribbon surface-link of total genus $g$, then $GS_0(F)=2^g$.

\phantom{x}

\noindent{\it Proof.} 
If $F$ is a ribbon surface-link of total genus $g$, then it is known 
that$GS_0(F)=2^g$, \cite{[8]}. 
For a self-trivial 1-handle $h_0$ on a surface-link $F$, 
the group $\Delta(F(h_0);Z_2)$ is the direct sum $\Delta(F;Z_2)\oplus Z_2\oplus Z_2$ 
or $\Delta(F;Z_2)\oplus Z_2$ according to whether a simple loop $\ell$ in $F$ induced from the core arc $\alpha_0$ of $h_0$ bounds an immersed disk $d$ in $S^4$ with 
$d\cap F=\ell$. In either case,  it is calculated that $GS_0(F(h_0))=GS_0(F)\cdot 2$. 
If $F$ is a quasi-ribbon surface-link of total genus $g$, then there is a self-trivial 
1-handle system $h$ such that $F(h)$ is a ribbon surface-link of total genus $g+s$ 
so that $GS_0(F(h)) =GS_0(F)\cdot 2^s=2^{g+s}$,   for the number $s$ of 1-handles in $h$.  Thus, $GS_0(F)=2^g$. This completes the proof of Lemma~3.1.

 \phantom{x} 
 
By using this invariant $GS_0(F)$, the proof of Theorem~1.2 is done as follows. 

\phantom{x} 

\noindent{\it 3.2: Proof of Theorem~1.2.}  
Let $k\cup k'$ be a non-split link in the interior of a 3-ball $B$ such that $k$ and $k'$ are trivial knots. For the boundary 2-sphere $S^B=\partial B$ and the disk $D^2$ with boundary circle $S^1$,  
let $K$ be the surface-link of  torus-components $T= k\times S^1$  
and $T'= k'\times S^1$ in the 4-sphere $S^4$ with  
$S^4=B\times S^1 \cup S\times D^2$, which is a ribbon surface-link in $S^4$, \cite{[2]}.  Then, $GS_0(K)=2^2$.  
Since $k$ and $k'$ are trivial knots in $B$, the torus-knots $T$ and $T'$ are trivial torus-knots in $S^4$  by construction.  
Since $k\cup k'$ is non-split in $B$, there is a simple loop $t(k)$ in $T$ coming from the longitude 
of $k$ in $B$ such that $t(k)$ does not bound any disk not meeting $T'$ in $S^4$,  
meaning that there is a simple loop $c$ in $T$ unique up to isotopies of $T$  
which bounds a disk $d$ in $S^4$ not meeting $T'$, where $c$ and $d$ are given by $c=\{p\}\times S^1$  and $d=a\times S^1\cup \{q\} \times D^2$ for a simple arc $a$ in $B$ joining a point $p$ of $k$ to a  point $q$ in $S$ with $a\cap(k\cup k')=\{p\}$ and $a\cap S=\{q\}$. 
Regard the 3-ball $B$ as the product $B=B_1\times [0,1]$ for a disk $B_1$.  
Let $\tau_1$ be a diffeomorphism of the solid torus $B_1\times S^1$ given by one full-twist rounding the meridian disk $B_1$ one time along the $S^1$-direction. Let $\tau=\tau_1\times 1$ be the product diffeomorphism of $(B_1\times S^1)\times[0,1]=B\times S^1$ for the identity map $1$ of $[0,1]$.  
Let $\tau_{\partial}$ be the diffeomorphism of the boundary $S^B\times S^1$ of 
$B\times S^1$  
obtained from $\tau$ by restricting to the boundary,  
and the 4-manifold $M$ obtained from $B\times S^1$ and $S^B\times D^2$ by pasting the boundaries  
$\partial (B\times S^1)=S^B\times S^1$ and $\partial (S^B\times D^2)=S^B\times S^1$ by the diffeomorphism  $\tau_{\partial}$. Since the diffeomorphism $\tau_{\partial}$ of $S^B\times S^1$ extends  
to the diffeomorphism $\tau$ of $B\times S^1$, the 4-manifold $M$ is diffeomorphic to $S^4$.  
Let $K_M=T_M\cup T'_M$ be the surface-link of torus components $T_M$ and  $T'_M$  in the 4-sphere $M$,  arising from $K=T\cup T'$ in  $B\times S^1$.  
There is a meridian-preserving isomorphism $\pi_1(S^4\setminus K, x)\to \pi_1(M\setminus K_M, x)$  by van Kampen theorem.  
Let $(S^4,K')=(M,K_M)$. 
To construct a desired pair  $(F,F')$ of ${\mathbf F}$-links $F$ and  
$F'$ from $(K, K')$, let $A$ be a 4-ball in $S^4$ such that 
$A\cap K=A\cap K'$ is a trivial disk pair 
in $A$ taking  one disk from each component of $K$ and one disk from each component of $K'$, where connected-summing of a trivial surface-knot and/or adding a trivial surface-link component as a split component are made. 
By van Kampen theorem, it is shown that the fundamental groups 
$\pi_1(S^4\setminus F,x)$ and $\pi_1(S^4\setminus F',x)$ are the same group up to meridian-preserving isomorphisms. 
By construction, $F$ is a ribbon surface-link. 
To show that $F'$ is not a quasi-ribbon surface-link, note that 
the loop $t(k)$ in $T_M$ does not bound any disk not meeting $T'_M$ in $M$,  
so that the loop $c$ in $T_M$ is a unique simple loop up to isotopies of $T_M$ which  
bounds a disk $d_M=a\times S^1\cup D^2_M$ in $M$ not meeting $T'_M$, where  
$D^2_M$ denotes a proper disk in $S^B\times D^2$ bounded by the loop  
$\tau_{\partial}(\{q\}\times S^1)$.  
An important observation is that the self-intersection number  
$\mbox{Int}(d_M, d_M)$ in $M$ with respect to the surface-framing on $K_M$ is $\pm1$. This means that the loop $c$ in $T_M$ is a non-spin loop.  
Similarly, there is a unique non-spin loop $c'$ in $T'_M$ which bounds  
a disk $d'_M$ with the self-intersection number  
$\mbox{Int}(d'_M, d'_M)=\pm 1$ with respect to the surface-framing on $K'=K_M$. 
Then it is calculated that $GS_0(K')=0$, so that $GS_0(F')=2^{g-2}$ for the total genus $g(\geq 2)$ of $F'$. Thus, $F'$ is not a quasi-ribbon surface-link by Lemma~3.1.  
This completes the proof of Theorem~1.2.  

\phantom{x}

In 3.2, the non-quasi-ribbon surface-link $K'$ of two components starting from 
the Hopf link 
$k \cup k'$ in the interior of a 3-ball $B$ has the free abelian fundamental group of rank $2$. Then by van Kampen theorem, the surface-knot $K'(h')$ obtained from $K'$ by surgery along any fusion 1-handle $h'$ on $K'$ has the infinite cyclic fundamental group, so that $K'(h')$ is a trivial surface-knot in $S^4$ by smooth unknotting result of a surface-knot, \cite{[3], [5]}. Since a non-quasi-ribbon surface-link is a non-ribbon surface-link, this shows that there is a non-ribbon surface-link $F'$ 
of ribbon components such that the 
surface-knot obtained from  $F'$  by surgery along any fusion 1-handle is a ribbon surface-knot,  meaning that {\it non-ribbonness of a surface-link cannot detect 
in general  by surgery along any fusion 1-handle system}, giving 
a strong counterexample to \cite[Theorem 1.4]{[4]}. The diffeomorphism 
$\tau_{\partial}$  of $S^B\times S^1$ in 3.2 coincides with Gluck's non-spin diffeomorphism of $S^2 \times S^1$, \cite{[13]}. The surface-link $(M, K_M)$ called a 
{\it turned torus-link} of a link $k \cup  k'$ in $B$ is an analogy of a 
{\it turned torus-knot} of a knot in $B$, \cite{[14]}. 

 \phantom{x} 

\noindent{\bf 4. Proof of Theorem~1.3} 

For a 1-handle $h$  on a surface-link $F$ in $S^4$, 
a trivial $S^2$-link {\it around} $h$ is 
a trivial $S^2$-link $O^h$ in $S^4$ all of whose components are  the boundary-spheres of fibers of  a trivial normal 3-disk bundle of the core arc $\alpha$ of $h$ in $S^4$  
not meeting with $F\cup h$. 
The following lemma concerns a replacement of a spanning arcof a surface-link.

\phantom{x} 

\noindent{\bf Lemma~4.1.}  Let $F$ be a surface-link  in $S^4$, and $\alpha_0$ a simple arc  spanning $F$ in $S^4$ meeting $F$ only with the endpoints. Then every  simple arc 
$\alpha$ spanning $F$ in $S^4$ meeting $F$ only with the same endpoints as 
$\alpha_0$ is isotopic to an arc obtained from the  arc $\alpha_0$ by band summing with  a meridian loop system $m$ of $F$ by a smooth isotopy of $S^4$  keeping $F$ fixed.

\phantom{x}

\noindent{\it Proof.}  
After a slight move of the interior of $\alpha$, the union $\alpha_0 \cup \alpha$ bounds a disk $D$ smoothly embedded in $S^4$. The interior of the disk $D$ transversely meets $F$ in a finite point set 
$Q=\{q_1, q_2, \dots, q_n\}$, whose regular neighborhood in $D$ is an unoriented meridian disk system $N(Q)=\{N(q_1), N(q_2), \dots, N(q_n)\}$ of  $F$ in $S^4$. The arc $\alpha$ is $\partial$-relatively isotopic  in $D$ to a band sum of the arc  $\alpha_0$ and the unoriented meridian loop system $m=\partial N(Q)=\{m_1, m_2,\dots,m_n \}$, where $m_i=\partial N(q_i)$, along a band system $b = \{b_1, b_2, \dots, b_n\}$ in $D$,  where $b_i$ spans 
$\alpha$ and $m_i$. This completes the proof of Lemma~4.1.  

\phantom{x} 

The following lemma is obtained from Lemma~4.1 by replacing the spanning arcs 
$\alpha_0$ and $\alpha$  with the 1-handles $h_0$ and $h$ on $F$ with core arcs 
$\alpha_0$ and $\alpha$, respectively.

\phantom{x}

\noindent{\bf Lemma~4.2.} 
Let $F$ be a surface-link  in $S^4$, and $h_0$ a  1-handle on $F$. Then for every 
1-handle $h$ on $F$ with the same attaching part as $h_0$, there is a trivial 
$S^2$-link  $O$ around $h_0$ in $S^4$ with an equivalence  
$F(h)=(F(h_0)\cup O)(\xi^O)$ 
for a fusion 1-handle system $\xi^O$ connecting the surface-link $F(h_0)$ and the trivial $S^2$-link $O$.

\phantom{x}

\noindent{\it Proof.}
In Lemma~4.1,  let $B = \{B_1, B_2, \dots, B_n\}$ be a disjoint 3-ball system  in $S^4$ attaching to $F$ with a disk system $d = \{d_1, d_2, \dots, d_n\}$ and meeting  a meridian loop system $m$ of $F$ transversely  so that $m \cap B = \{t_1, t_2, \dots, t_n\}$ for a point $t_i = m_i \cap (B_i\setminus\partial B_i)\, (i=1,2,\dots,n)$. 
The thickening 1-handle $h$ of $\alpha$ meets $B_i$ with a transverse disk 
centered at $t_i$. 
For the complementary disk system $d^c = \mbox{cl} (\partial B \setminus d)$ of 
$d$ in the boundary $\partial B$ of $B$, regard the surface-link 
$F'= \mbox{cl}(F\setminus d)\cup d^c$ as $F$.
Deform the arc $\alpha$ into the arc $\alpha_0$  by taking together the 3-ball system $B$ along the band system $b$. 
Then the 1-handle $h$ is deformed into the 1-handle $h_0$ thickening the arc $\alpha_0$ by taking together the 3-ball system $B$. This means that there is an equivalence 
$F(h)=(F(h_0)\cup O)(\xi^O)$ for the trivial $S^2$-link  $O=\partial B$ around $h_0$ and  a fusion 1-handle system  $\xi^O$ connecting  $F(h_0)$ and $O$ 
which is constructed by thickening and stretching the complementary disk system $d^c$ along the band system $b$.  This completes the proof of Lemma~4.2.

\phantom{x}

For a surface-knot $F$ in $S^4$, let $V$ be a connected Seifert hypersurface of $F$ in $S^4$, and $(V,F)\times[0,1]$ a normal $[0,1]$-bundle over $(V,F)$ in $S^4$.
Let $(V_i,F_i) =(V,F)\times \{i\}\, (i=0,1)$. 
Then the natural homomorphism $H_1(F_1; Z) \to H_1(S^4\setminus F_0; Z)$ is the zero map.  The surface-link $P(F)=F_0\cup F_1$  is called the {\it anti-parallel surface-link} of $F$, where, by convention, $F_0$ and $F_1$ are identified with  
$-F$ (i.e., the orientation reversed $F$) and $F$, respectively. 
Cochran's conjecture says that when $F$ is a non-ribbon $S^2$-knot, 
the anti-parallel $S^2$-knot $P(F;h)$ must be a non-ribbon $S^2$-knot 
for a sufficiently complicated fusion 1-handle $h$, \cite{[10]}. 
Let $h$ be a fusion 1-handle on $P(F)$, and 
$\beta$  a 1-handle system on $F$ embedded in $V$ 
 with $V(\beta)=\mbox{cl}(V\setminus \beta)$ a handlebody such that 
$\beta_i=\beta\times\{i\}$ is embedded in $V_i\, (i=0,1)$ and disjoint from $h$. 
A fusion 1-handle $h$ on $P(F)$ is {\it surgically trivial}, or simply {\it S-trivial} if 
$P(F)(\beta_0\cup\beta_1)(h)$ is a trivial surface-knot. 
Otherwise, $h$ {\it surgically nontrivial}, or simply {\it S-nontrivial}.
Since $P(F)$ is a boundary surface-link, this definition does not depend on 
choices of $V$ and $\beta$, \cite[Lemma~3.1]{[11]}. 
Cochran's conjecture says that  for the anti-parallel $S^2$-link $P(K)$ of a non-ribbon 
$S^2$-knot $K$ in $S^4$, if $h$ is a \lq\lq sufficiently complicated\rq\rq fusion 1-handle 
$h$ on $P(K)$ must be a non-ribbon$S^2$-knot, \cite{[10]}. 
By replacing \lq\lq sufficiently complicated\rq\rq with \lq\lq S-nontrivial\rq\rq ,
the following theorem not only confirms Cochran's conjecture affirmatively but also 
shows that the result holds even when the non-ribbon $S^2$-knot $K$ is generalized to 
any non-ribbon surface-knot $F$, which was announced in \cite{[15]} 
with incomplete proof.

\phantom{x} 

\noindent{\bf Theorem 4.3.} Let $P(F)$ be the anti-parallel surface-link 
of a non-ribbon surface-knot $F$ in $S^4$. 
According to whether a fusion 1-handle $h$ on $P(F)$ is S-trivial or S-nontrivial, the surface-knot $P(F)(h)$ is a trivial surface-knot or a non-ribbon surface-knot, respectively.  

\phantom{x} 

\noindent{\it Proof.} 
 For any S-nontrivial fusion 1-handle $h$ on $P(F)$, the surface-knot $P(F)(h)$ is a 
non-ribbon surface-knot,  \cite{[11]}.  
Let $h$ be  an S-trivial fusion 1-handle  on $P(F)$.
Let $V$ be a connected Seifert hypersurface for $F$, and 
$\beta$  a 1-handle system on $F$ embedded in $V$ such that 
$\mbox{cl}(V\setminus\beta)$ is a handlebody. 
Let  $(V,F)\times [0,1]$ be a normal $[0,1]$-bundle over $(V,F)$ in $S^4$.
For a subset $X$ in $V$, let $X_i$ denote $X\times\{i\}$ for $i=0$ of $1$. 
Let $P(F)=F_0\cup F_1$. 
For a disk $d$ in $F$, let
$h_0 = d \times[0,1]$ be the standard trivial 1-handle on $P(F)$, and assume that the attaching part of $h$ to $P(F)$ coincides with the disk union $d_0\cup d_1$. 
Note that the surface-knot $P(F)(h_0)$ is a trivial surface-knot since it bounds a handlebody $\mbox{cl}(F\setminus d)\times[0,1]$. 
Let $\alpha$ and $\alpha_0$ be the core arcs of the 1-handles
$h$ and $h_0$ on $P(F)$, respectively, which are disjoint except for the same endpoints $p_i$ in $d_i\,(i=0,1)$.  The simple loop $\alpha_0\cup \alpha$  meets 
$V_0\cup V_1$ only with the points $p_i\, (i = 0, 1)$, and bounds 
a smoothly embedded disk $\Delta$ in $S^4$, \cite[Lemma~3.1]{[11]}, \cite{[12]}.
Let $D^{\beta}$  be a transverse disk system  of 
the 1-handle system $\beta$ with just one disk for every 1-handle. 
By regarding the handlebody systems $V_i(\beta_i)\, (i=0,1)$ as spine  graphs 
attaching to the points $p_i\, (i = 0, 1)$, the disk $\Delta$ can avoid intersecting 
the handlebody system $V_0(\beta_0)\cup V_1(\beta_1)$ except for the points 
$p_i\, (i = 0, 1)$ by general position. In particular, the disk $\Delta$ does not meet the trivial surface-link $P(F)(\beta_0\cup \beta_1)$ except for the points 
$p_i\, (i = 0, 1)$. Thus, the disk $\Delta$ meets transversely the disk system 
$D_0^{\beta}\cup D_1^{\beta}$ with a finite number of  interior points.  
By  Lemma~4.1 applied to the disk $\Delta$, 
the arc $\alpha$ is obtained from the arc $\alpha_0$ by band-summing a meridian loop system $m$  of  the disk system $D_0^{\beta}\cup D_1^{\beta}$ in $S^4$. 
By  Lemma~4.2,  there is a disk system
$E_i^{\beta}$ which is obtained from $D_i^{\beta}$  and 
 a trivial $S^2$-link $O_i$  around the 1-handle $h_0$ by surgery 
along a fusion 1-handle system $\xi_i$ such that 
the surface-knot 
$P(F)(h)(\beta_0\cup \beta_1)(D_0^{\beta}\times I\cup D_1^{\beta}\times I)$ is equivalent to the surface-knot 
$P(F)(h_0)(\beta_0\cup \beta_1)(E_0^{\beta}\times I\cup E_1^{\beta}\times I)$. 
By choosing a disjoint proper arc system $\gamma$ in $\mbox{cl}(F(\beta)\setminus d)$ 
such that each arc of $\gamma$ meets transversely the loop system 
$\partial D^{\beta}$ with just one point, the disk system 
$\Lambda=\gamma\times[0,1]$ is considered as a core disk system 
of a 2-handle $\Lambda\times I$ on $P(F)(h_0)(\beta_0\cup \beta_1)$ 
such that  the pairs $(E_i^{\beta}\times  I, \Lambda\times I)\,(i=0,1)$ and 
$(D_i^{\beta}\times  I, \Lambda\times I)$ are O2-handle pair systems on 
$P(F)(h_0)\beta_0\cup \beta_1)\, (i=0,1)$ with the same attaching parts for each $i$, where the 1-handle systems $\xi_i\,(i=0,1)$ are chosen to be disjoint from the disk 
$\Lambda$. 
By common 2-handle property, the O2-handle pair system 
$(E_i^{\beta}\times  I, \Lambda\times I)$ is sent to  the O2-handle pair system 
$(D_i^{\beta}\times  I, \Lambda\times I)$  for each $i$ 
by an orientation-preserving diffeomorphism $f_i$ of $S^4$ 
keeping $P(F)(h_0)(\beta_0\cup \beta_1)$ fixed, \cite{[3],[5]}. 
This implies that the surface-knot $P(F)(h)=
P(F)(h_0)(\beta_0\cup \beta_1)(E_0^{\beta}\times  I\cup E_1^{\beta}\times  I)$
is equivalent to 
$P(F)(h_0)(\beta_0\cup \beta_1)(D_0^{\beta}\times  I\cup D_1^{\beta}\times  I)
=P(F)(h_0)$, which is a trivial surface-knot.  
This completes the proof of Theorem~1.2. 

\phantom{x} 

The proof of Theorem~1.3 is done as follows. 

\phantom{x} 

\noindent{\it  4.4: Proof of Theorem~1.3.} 
Let $K$ be a non-ribbon $S^2$-knot in $S^4$, and  $V$ a connected Seifert hypersurface 
for $K$. Let $B\cup V^B$ be a decomposition of $V$ with $B$ a 3-ball 
and $d^B=B\cap V^B$ a proper disk in $V$.  For a collar $V\times [0,1]$ of $V$ in $S^4$, 
let  $h_0=d\times[0,1]$ be the trivial fusion 1-handle on   
the anti-parallel $S^2$-link $P(K)=K_0\cup K_1$ of $K$ for a disk $d$  
in the disk $\mbox{cl}(\partial B\setminus d^B)$, and $h$ 
a non-trivial fusion 1-handle on $P(K)$ with the same attachment as $h_0$ 
such  that $h$ is in a  regular neighborhood $N$ of the 3-ball 
$B_0\cup h_0\cup B_1$ in $S^4$ and meets the interior of the 3-ball 
$B_0$ with $n$ transversal disks of $h$, for any given positive integer $n$. 
Let $h'$ be a fusion 1-handle on $P(K)$ in $N$ with the same attachment as $h$ 
such that the core arc $\alpha'$ of $h'$ does not meet $B_0$ except for the attaching point and the core arc $\alpha$ of $h$ is a band sum of  $\alpha'$ 
and a meridian loop system on $K_0$ by an argument of Lemma~4.1. 
Since $\alpha'$ is $\partial$-relatively isotopic to $\alpha_0$ in $N$, the $S^2$-knot 
component $P(K)(h')$ is equivalent to  the trivial $S^2$-knot $(P(K)(h_0)$. 
By Lemma~4.2,  there is an equivalence  
$P(K)(h)=(P(K)(h_0)\cup O)(\xi^O)$ 
for a trivial $S^2$-link $O$ with $n$-components around $h_0$  and a fusion 1-handle system $\xi^O$ connecting the  surface-knot $P(K)(h_0)$ and $O$.  
By Theorem~4.3,  the $S^2$-knot $P(K)(h)$ is non-ribbon, so that 
the $S^2$-link $F_K=P(K)(h_0)\cup O$ is a non-ribbon $S^2$-link of $n+1$ trivial components. By Theorem~1.1,  $F_K$ is quasi-ribbon.  
For every disconnected closed oriented surface ${\mathbf F}$, a quasi-ribbon 
but non-ribbon ${\mathbf F}$-link $F$ with trivial components 
is constructed  from $F_K$ for some $n$ 
by connected-summing trivial surface-knots  and/or adding a trivial surface-link component  as a splitting component. Non-ribbonness of $F$ is seen from that if $F$ is a ribbon surface-link, then $F_K$ must be a ribbon $S^2$-link, \cite{[4]}.  
To see that the ${\mathbf F}$-link $F$ is quasi-ribbon, 
note that $F_K$ is quasi-ribbon and by Theorem~2.1 there is a self-trivial 
1-handle system $h^*$ on $F_K$ such that $F_K(h^*)$ is a ribbon surface-link.
The connected-summands of trivial surface-knots  and/or addition of 
 a trivial surface-link component are taken locally apart from $h^*$, and 
 $F(h^*)$ is taken to be a ribbon surface-link. By Theorem~2.1, $F$ is 
 a quasi-ribbon but non-ribbon ${\mathbf F}$-link of trivial components. 
  This completes the proof of Theorem~1.3.

\phantom{x}

\noindent{\bf Conclusion}  

The original question is when a surface-link $F$ of ribbon surface-knot components becomes a ribbon surface-link, \cite{[4]}. 
It is shown that a boundary surface-link $F$ is a ribbon surface-link if  
there is a local S-nontrivial fusion 1-handle system $h$ on $F$ with $F(h)$ a ribbon surface-link, \cite[Theorem~1.1]{[11]}.  For surface-links with at least two non-sphere components, there are 
examples of non-ribbon surface-links $F$ of ribbon components. In fact, there are 
examples of  non-ribbon surface-links of trivial components which 
cannot be detected by the fundamental groups, as shown in Theorem~1.2 as well as   \cite{[6]}. According to Theorem~1.1, a surface-link of trivial components with at most one 
non-sphere component is a quasi-ribbon surface-link, which shares properties analogous to those of ribbon surface-links. However, as Theorem~1.3 shows, additional conditions are required for it to be a ribbon surface-link.

\phantom{x} 

\noindent{\bf Acknowledgements}  

This work was partly supported by JSPS KAKENHI Grant Numbers JP21H00978 and JP26K06456, and MEXT Promotion of Distinctive Joint Research Center Program JPMXP0723833165 and Osaka Metropolitan University Strategic Research Promotion Project (Development of International Research Hubs).


\begin{thebibliography}{99} 

\bibitem{[1]} Kawauchi, A., Shibuya, T., Suzuki, S. (1982). Descriptions on surfaces in four-space,  II: Singularities and cross-sectional links, Math Sem Notes Kobe Univ, 11: 31-69. 

\bibitem{[2]} Kawauchi, A. (2015). A chord diagram of a ribbon surface-link, J Knot Theory Ramifications, 24: 1540002 (24 pages). 

\bibitem{[3]} Kawauchi, A. (2021). Ribbonness of a stable-ribbon surface-link, I. A stably trivial surface-link, Topology and its Applications, 301:  107522 (16pages).

\bibitem{[4]} Kawauchi, A. (2025). Ribbonness of a stable-ribbon surface-link, II:  General case, (MDPI) Mathematics, 13 (3): 402 (1-11).

\bibitem{[5]} Kawauchi, A. (2023).  Uniqueness of an orthogonal 2-handle pair on a surface-link, Contemporary Mathematics (UWP), 4: 182-188.

\bibitem{[6]} Kawauchi, A. (2024). Note on surface-link of trivial components, 
J Comp Pure Appl Math, 2 (1): 1 - 5.

\bibitem{[7]} Kawauchi, A. (2025).  Revised note on surface-link of trivial components.  J Comp Pure Appl Math, 3 (2): 1-10. 
 
\bibitem{[8]} Kawauchi, A. (2002). On pseudo-ribbon surface-links, J Knot Theory Ramifications, 11: 1043-1062.

\bibitem{[9]} Ogasa, E. (2001). Nonribbon 2-links all of whose components are trivial knots and some of whose band-sums are nonribbon knots, J Knot Theory Ramifications 10: 913-922.

\bibitem{[10]} Cochran, T. (1983). Ribbon knots in $S^4$, J London Math Soc (2) 28: 563-576.  

\bibitem{[11]} Kawauchi, A. (2026). Fusion of boundary surface-link and ribbonness,
  J Comp Pure Appl Math, 4 (1): 1-13. 

\bibitem{[12]} Hosokawa, F. and  Kawauchi, A. (1979). Proposals for unknotted 
surfaces in four-space, Osaka J Math, 16: 233-248. 

\bibitem{[13]} Gluck, H. (1962). The embedding of two-spheres in the four-sphere, Trans Amer Math Soc, 104: 308-333. 

\bibitem{[14]}  Boyle, J. (1993). The turned torus knot in $S^4$, J Knot Theory Ramifications, 2: 239-249.

\bibitem{[15]} Kawauchi, A. (2025). Ribbonness on boundary surface-link. Open Access J Phys Math, 1 (2): 1-4.   


 
\end{thebibliography}
\end{document}